\documentclass[12pt]{article}
\usepackage{amsmath2000,amsthm2000,amssymb,amsfonts}

\makeatletter
\newfont{\footsc}{cmcsc10 at 8truept}
\newfont{\footbf}{cmbx10 at 8truept}
\newfont{\footrm}{cmr10 at 10truept}
\renewcommand{\ps@plain}{%
\renewcommand{\@oddfoot}{\footsc the electronic journal of combinatorics
  {\footbf 9} (2002), \#R10\hfil\footrm\thepage}}
\renewcommand{\section}{\@startsection
{section}{1}{0pt}{-3.5ex plus-1ex minus-.2ex}%
{0.8ex plus.2ex}{\centering\normalfont\large\bfseries}}
\newcommand{\ssection}[1]%
 {\addtocounter{section}{1}\section*{\arabic{section}.\ \ #1} \ }
\newcommand{\ssections}[1]%
 {\addtocounter{section}{1}\section*{#1} \ }
\makeatother
\theoremstyle{plain}
\newtheorem{thm}{Theorem}
\newtheorem{cor}[thm]{Corollary}
\newtheorem{lem}[thm]{Lemma}
\newtheorem{pr}[thm]{Proposition}

\numberwithin{thm}{section}

\numberwithin{equation}{section}

\newcommand{\Z}{{\mathbb Z}}
\newcommand{\N}{{\mathbb N}}
\newcommand{\E}{{\bf E}}

\newcommand{\deq}{\, {\stackrel {\textup{def}} {=}}}

\newcommand{\ee}{\varepsilon}

\newcommand{\cS}{{\mathcal S}}

\newcommand{\Prob}{{\bf{P}}}
\newcommand{\F}{{\mathcal{F}}}

\title{Generating a random sink-free orientation\\
in quadratic time}

\author{\begin{tabular}{ccc}
Henry Cohn
\!\!\phantom{\thanks{This work was begun at the 1997 Workshop on
Exact Simulation in Rebild, Denmark (funded by the European
Science Foundation) and completed at the 1998 Institute for
Elementary Studies (funded by the National Science Foundation).
Cohn was supported by an NSF Graduate Research Fellowship, and
currently holds a five-year fellowship {}from the American
Institute of Mathematics. Propp is supported in part by grants
{}from the National Science Foundation and National Security
Agency.  Pemantle is supported in part by NSF grant DMS 9803249.}}
& Robin Pemantle & James Propp\\
\small Microsoft Research &
\small Department of Mathematics & \small
Department of Mathematics\\[-0.8ex]
\small One Microsoft Way & \small Ohio State University
& \small University of Wisconsin\\[-0.8ex]
\small Redmond, WA 98052-6399 & \small Columbus, OH 43210
& \small Madison, WI 53706\\[-0.8ex]
\small \texttt{cohn@microsoft.com} & \small
\texttt{pemantle@math.ohio-state.edu}
& \small \texttt{propp@math.wisc.edu}\\[-0.8ex]
\end{tabular}}

\date{\small Submitted: May 1, 2001;  Accepted: March 12, 2002.
\\ \small MR Subject Classifications: 60C05, 68W20}

\begin{document}

\maketitle

\begin{abstract}
A {\em sink-free orientation\/} of a finite undirected graph is a
choice of orientation for each edge such that every vertex has
out-degree at least 1.  Bubley and Dyer (1997) use Markov Chain
Monte Carlo to sample approximately {}from the uniform
distribution on sink-free orientations in time $O(m^3 \log (1 /
\ee))$, where $m$ is the number of edges and $\ee$ the degree of
approximation. Huber (1998) uses coupling {}from the past to
obtain an exact sample in time $O(m^4)$.  We present a simple
randomized algorithm inspired by Wilson's cycle popping method
which obtains an exact sample in mean time at most $O(nm)$, where
$n$ is the number of vertices.
\end{abstract}

\ssection{Introduction}

A common problem is to select a random sample efficiently {}from a
large collection of combinatorial objects.  There are many
reasons one may wish to do this. One is to obtain an approximate
count: Jerrum and Sinclair \cite{JS} showed that if one can
generate nearly uniform samples, then for each $\ee >0$, one can
obtain the cardinality of the collection to within a factor of $1
+ \ee$ with probability $1 - \ee$, in just a little more time.
When counting the collection is \#P-hard, as in the case of
properly $k$-coloring a graph, this may be the only reasonable
way to count, since it is unlikely that \#P-hard counting
problems can be solved exactly in polynomial time. Another reason
to seek a sampling algorithm is that it may shed light on
properties of the typical sample. For example, the analysis of
typical spanning trees of Cayley graphs \cite{P,BLPS} relies on
two algorithms, the first developed by Aldous \cite{A} and Broder
\cite{B} and the second by Wilson \cite{W}.  The analysis of
phase boundaries in typical domino tilings of regions known as
Aztec diamonds also relies on a sampling algorithm, known as
domino shuffling \cite{CEP}. Finally, sample generation may be a
way of producing conjectures about the typical sample via
simulation, when no theorem is known (for example, the results in
\cite{CEP} were initially discovered this way).

One common way to generate samples is Markov Chain Monte Carlo
(MCMC).  Here one finds an ergodic Markov chain whose equilibrium
measure is the desired distribution $\mu$; then one runs the
chain until the distribution is close to $\mu$.  Constructing
such a chain is usually easy (often when $\mu$ is uniform, there
is a natural doubly stochastic transition matrix) and the hard
part is knowing how long to run it.  This may be established via
eigenvalue bounds, or via coupling arguments or stopping times.
In cases where the time bounds on the chain are established via
coupling, it is often possible to improve on MCMC by using
coupling {}from the past (CFTP) to obtain an exact sample rather
than an approximate one \cite{PW}.

In this note we consider the generation of a random {\em sink-free
orientation} (SFO) of a finite undirected graph.  Sink-free
orientations were introduced by Bubley and Dyer \cite{BD}, who
were motivated by an equivalence between counting them and
counting satisfying assignments of Boolean formulas in
conjunctive normal form in which each variable occurs at most
twice (they call this problem Twice-SAT). Bubley and Dyer showed
that counting sink-free orientations is \#P-complete, so it is
unlikely that an exact count can be obtained in polynomial time,
and we must use approximate counting techniques based on nearly
uniform sampling.

Bubley and Dyer give an MCMC algorithm that produces a sample
whose distribution is within $\ee$ of uniform (in total variation)
in time $O(m^3 \log (1/\ee))$, where $m$ is the number of edges.
Huber \cite{H} uses Bubley and Dyer's analysis along with CFTP to
produce an exact uniform sample in mean time $O(m^4)$. The
purpose of this note is to improve the running time to $O(nm)$,
where $n$ is the number of vertices.  Instead of MCMC, we use a
strong uniform time algorithm inspired by David Wilson's cycle
popping algorithm \cite{W} for generating uniform directed
spanning trees.

We now describe the problem and our results more precisely. Let
$G = (V, E)$ be a finite undirected graph. We allow multiple edges
and self-loops (but at most one self-loop per vertex, since
multiple self-loops play no useful role in sink-free
orientations). We define an $n$-cycle to be a ring of $n$ vertices
$v_0,\dots,v_{n-1}$ with edges {}from $v_i$ to $v_{i+1}$ for each
$i$ (taken modulo $n$;  note that a $1$-cycle is a vertex with a
self-loop), and an $n$-lollipop to be a path consisting of $n$
vertices and $n-1$ edges, with a self-loop added at one end.

An orientation of an edge between vertices $v$ and $w$ is a
mapping of the set $\{\textup{head},\textup{tail}\}$ onto $\{ v ,
w \}$.   Thus, a self-loop has only one orientation, but all
other edges have two.  To reverse the orientation of an edge,
swap its head and tail. An orientation of $G$ is an orientation of
each edge. A sink in an orientation is a vertex that is not the
tail of any edge (a source is the opposite, i.e., not the head of
any edge), and a sink-free orientation (SFO) of $G$ is an
orientation that contains no sinks. If any connected component of
$G$ is a tree, then $G$ has no SFO, and vice versa. Henceforth we
restrict consideration to the class $\cS$ of graphs in which no
component is a tree.  Let $\mu_G$ denote the probability measure
assigning probability $1/N$ to each SFO of $G$, where $N$ is the
total number of SFO's of $G$.

Our algorithm, which we call ``sink popping,'' works as follows.
Given a graph, orient the edges by independent, fair coin flips.
If this orientation has no sinks, then it is the SFO we seek.
Otherwise, choose any sink, and randomly re-orient each edge that
points into the sink (i.e., all of its edges). We call this
popping the sink, for reasons that will become clearer in the
next section. Repeat until there are no more sinks.

We now state our main result.

\begin{thm} \label{th:main}
For every graph $G \in \cS$, sink popping terminates in finite
time with probability~$1$, regardless of how one chooses which
sink to pop, and produces an output whose distribution is
precisely $\mu_G$. The average number of sinks that must be
popped is at most $\binom{n}{2}$, where $n$ is the number of
vertices of $G$, once again regardless of how one chooses which
sink to pop.  Equality holds only for the $n$-cycle or the
$n$-lollipop. The expected number of times each particular vertex
is popped is at most $n-1$.
\end{thm}

Sink popping is briefly mentioned at the end of \cite{PW2}, where
the claims in the first sentence of the theorem are mentioned
without detailed proof (and the running time is not analyzed).

To show that sink popping's running time is $O(nm)$, we need to
state the algorithm slightly more carefully.  The subtle point is
avoiding spending lots of time searching for sinks.  We will keep
a list of all sinks in the graph, and also a table showing the
out-degree of every vertex. Generating these initially takes time
proportional to the sum of the degrees of the vertices, or $O(m)$
time.  Whenever we search for a sink to pop, we simply take the
first sink {}from the list. When we pop the sink, we update the
table to reflect the changes to its out-degree, and to those of
its neighbors.  It neighbors may have become sinks, in which case
we append them to the list of sinks.  (The purpose of the table
is to let us easily see whether the neighbors have become sinks,
without having to examine all their edges: in a complete graph,
that would waste lots of time.) No sink can be annihilated except
the one we popped, since no two sinks can share a common edge (it
would have to point to both). Thus, each time we pop a sink at
$v$, re-orienting its edges and updating the list and table
requires time $O(\deg(v))$.  By Theorem~\ref{th:main}, the
expected number of times $v$ is popped is $O(n)$, so the total
expected number of operations is
$$
O\left(\sum_{v} n\deg(v) \right) = O(nm).
$$

This time bound does not actually estimate the number of bit
operations, but instead treats individual graph operations as
units.

In the next section we give another description of the sink
popping algorithm and explain its connection to cycle popping.
We also state some further results about sink popping with
arbitrary initial conditions.  The third section contains proofs
of the diamond and strong uniformity lemmas, which are analogous
to the equivalent lemmas for cycle popping. The fourth section
analyzes the running time.  The fifth section derives some
further facts about the running time.  We conclude with some
speculations and open questions.

\ssection{Sink popping and cycle popping}

Let $H = (V, E)$ be a finite, connected, directed graph, and let
$v$ be any fixed vertex of $H$.  A directed spanning tree of
$(H,v)$ is a subset of edges so that every vertex other than $v$
has out-degree 1 and $v$ has out-degree 0.  Wilson \cite{W,PW2}
invented the following algorithm, known as ``cycle popping,'' for
generating a uniform random directed spanning tree.  For each $w
\in V \setminus \{ v \}$ and $k \geq 0$, let $X_{w , k}$ be a
random edge leading out of $w$, chosen uniformly {}from among all
edges leading out of $w$.  Let these be independent as $w$ and
$k$ vary.  For fixed $w$, imagine the collection $\{ X_{w , k} : k
\geq 0 \}$ as a stack with $X_{w,0}$ on top.  Initially, look at
the collection $\{ X_{w , 0} : w \in V \setminus \{ v \} \}$,
that is, consider the collection $\{ X_{w , f(w)} \}$ with $f
\equiv 0$. If these form a directed spanning tree, stop and set
the sample equal to it.  If not, there must be a cycle in this
collection. Choose a cycle (it doesn't matter which), and
increment $f(w)$ by 1 for each $w$ in the cycle. (Imagine popping
these edges off the stack so the next element of each stack is now
on top.)  If the collection $\{ X_{w , f(w)} : w \in V \setminus
\{ v \} \}$ is now a directed spanning tree, stop and return this
for your sample, otherwise continue popping until you do stop.
Wilson showed that the set of cycles popped does not depend on
which you choose to pop when you have a choice, and that the
algorithm stops almost surely at a directed spanning tree with
uniform distribution.

We can describe sink popping in similar terms, which will be
useful in the proof of Theorem~\ref{th:main}. Let $G = (V , E)$
be a finite undirected graph in $\cS$ and let $\Omega =
\Omega_0^{\N}$, where $\Omega_0$ is the set of orientations of
$G$, i.e., $\Omega$ consists of sequences of orientations of $G$.
We endow $\Omega$ with the $\sigma$-field $\F$ generated by the
coordinate functions $X_{e , k}$ for $e \in E(G)$ and $k \geq 0$,
which specify the orientation of $e$ in the $k$-th orientation in
the sequence. Endow $\Omega$ with the probability measure $\Prob$
under which the coordinate functions are independent and each
equally likely to yield either orientation.  The intuition is
that $\{ X_{e,k} : k \geq 0 \}$ represents a stack of arrows
under the edge $e$. Define a random function $f : E \times \N \to
\N$ as follows. Let $f(e,0) = 0$ for all $e$.  Given $f (e , k)$
for all $e$, define $f(\cdot , k+1)$ inductively:  If the
collection $\{ X_{e , f(e , k)} : e \in E \}$ is an SFO, then set
$f(e,k+1) = f(e,k)$ for all $e$.  If not, choose a sink $v_k \in
V$ arbitrarily, i.e., a vertex $v_k$ for which all edges $e$
incident to it are oriented toward it by the orientation $X_{e ,
f(e , k)}$.  Let $f(e,k+1) = f(e,k)$ for $e$ not incident to
$v_k$, and $f(e,k+1) = f(e,k) + 1$ for $e$ incident to $v_k$. The
dependence of $f$ on the choice rule (for choosing $v_k$, if
there are several sinks) is suppressed in the notation, as is the
dependence on the choice of $\omega \in \Omega$ via the variables
$\{ X_{e,k}\} $. Intuitively, $f(e,k)$ is the original depth of
the arrow under $e$ now at the top of the stack at time $k$.  Say
that $v_k$ is the sink popped at time $k$, and let $\tau = \min
\{ k : \{ X_{e , f(e,k)} \} \textup { is an SFO} \}$ be the number
of pops before an SFO is obtained (conceivably $\tau=\infty$),
and $\eta = \eta (\omega, \textup{choice rule})$ denote the
resulting SFO (if any).

Except for its last sentence, Theorem~\ref{th:main} is
established by showing that $\tau < \infty$ with probability~$1$,
the law of $\{ X_{e , f(e,\tau)} \}$ is precisely $\mu_G$, and
$\E \tau \leq \binom{n}{2}$, with equality in and only in the
cases indicated. The first two lemmas are analogous to those used
by Wilson \cite{W} in establishing the validity of the cycle
popping algorithm, and the third is the running time analysis.

\begin{lem}[Diamond lemma] \label{lem:diamond}
The number of pops $\tau \leq \infty$ in a maximal popping
sequence is independent of the choice rule, as is the multiset
$\{ v_k : 0 \leq k < \tau \}$.  If $\tau < \infty$ then the
resulting SFO $\eta$ is also independent of the choice rule.
\end{lem}

The name ``diamond'' is meant to remind the reader that moving
from the top of a diamond to the bottom by going southeast then
southwest is equivalent to going southwest then southeast. This
terminology comes from the article \cite{E}.

\begin{lem}[Strong uniform time] \label{lem:unif}
Let $N$ be the number of SFO's of $G$.  Then for each $k \geq 0$,
and each SFO $\eta$,
$$
\Prob (\tau = k , \{ X_{e , f(e ,\tau)} \} = \eta) =
\frac {\Prob (\tau = k)}   {N} \, .
$$
In other words, $\tau$ is a strong uniform time.
\end{lem}

\begin{lem} \label{lem:n choose 2}
If $G \in \cS$ has $n$ vertices, then $\E \tau \leq \binom{n}{2}$,
with equality only for the $n$-cycle and the $n$-lollipop.
\end{lem}

We conclude this section by stating two results that shed further
light on the running time of the popping algorithm.

\begin{pr} \label{pr:equal}
The distribution of $\tau$ for the $n$-cycle is exactly the same
as the distribution for the $n$-lollipop.
\end{pr}

\begin{pr} \label{pr:extremal}
Let $G \in \cS$ be any graph with $n$ vertices and let $\F_0$ be
the $\sigma$-field generated by the variables $\{ X_{e , 0} \}$.
Then the conditional mean running time $\E (\tau \,|\, \F_0)$ is
always bounded by $n (n-1)$, and the only case to achieve this is
an $n$-lollipop with all edges oriented opposite to their
orientation in the unique SFO.
\end{pr}

\ssection{Strong uniformity}

We first establish deterministic facts holding for every sample
$\omega \in \Omega$. Say that a sequence $v_0 , \ldots , v_{k-1}$
with $k \leq \infty$ is a maximal popping sequence for $\omega$
if it is legal (i.e., only sinks are popped) and cannot be
extended to larger $k$ (thus if $k < \infty$ it results in an
SFO).  Note that if $k = \infty$, we do not mean our notation to
suggest that $v_0,v_1,\dots$ is followed by a final term
$v_{\infty-1}$;  instead, $v_0,\dots,v_{\infty-1}$ denotes the
infinite sequence $v_0,v_1,\dots$, with no final term.

Let $f(e,k)$ denote the function $f(e,k,\omega,\mathbf{v})$ where
$\omega$ is a sample point and $\mathbf{v}$ is a specified legal
sequence of pops of length at least $k$.  Define an equivalence
relation on finite sequences $v_0 , \ldots , v_{k-1}$ of vertices
of $G$ by calling two sequences equivalent if one can be changed
to the other by a sequence of transpositions of pairs of vertices
$(v_i , v_{i+1})$ that are not neighbors in $G$.  (Note that such
a transposition does not change whether a sequence is a legal
popping sequence.) The following lemma is useful, though obvious.

\begin{lem}[Deterministic strong Markov property] \label{DSMP}
Given an integer $j$ and vertices $v_0 , \ldots , v_{j-1}$, let
$\omega$ be any initial configuration for which $v_0 , \dots ,
v_{j-1}$ is a legal popping sequence and let $\omega'$ and
$\omega$ be related by
$$X_{e,k} (\omega') = X_{e , k + f(e,j)} (\omega) .$$
That is, $\omega'$ looks like $\omega$ after $v_0 , \ldots ,
v_{j-1}$ are popped.  Then the following deterministic strong
Markov property (DSMP) holds.  For any $k \le \infty$, the set of
sequences $\{ v_{j+i} : 0 \leq i \leq k-1 \}$ for which $v_0 ,
\ldots , v_{j+k-1}$ is a legal popping sequence for $\omega$ is
the same as the set of legal popping sequences of length $k$ for
$\omega'$.  If $v_0 , \ldots , v_{j+k-1}$ is maximal for $\omega$
then $v_j , \ldots , v_{j+k-1}$ is maximal for $\omega'$, and
leaves the same SFO (if $k < \infty$).
\end{lem}

Extend the definition of equivalence to infinite sequences by
saying that $v_0,v_1,\dots$ is equivalent to $w_0,w_1,\dots$ if,
by a sequence of transpositions as above applied to
$v_0,v_1,\dots$, one can transform $v_0,v_1,\dots$ so that
arbitrarily long initial segments of it match those of
$w_0,w_1,\dots$.  In particular, this implies that the multisets
$\{ v_k \}$ and $\{ w_k \}$ are the same.

Let $l(\omega)$ denote the minimal length of a maximal popping
sequence for $\omega$.

\begin{lem} \label{lem:equiv}
The set of maximal popping sequences is an equivalence class.
\end{lem}

\begin{proof}
Let $v_0 , \ldots , v_{k-1}$ be a legal popping sequence for
$\omega$ with $k < \infty$, and let $w_0 , \ldots , w_{k-1}$ be
obtained {}from $v_0 , \ldots , v_{k-1}$ by transposing $v_i$ and
$v_{i+1}$ which are not neighbors in $G$. Suppose $i=0$.  Since
the edges incident to $v_0$ are disjoint {}from the edges incident
to $v_1$, we may apply the DSMP to $v_0 , v_1$ and to $v_1 , v_0$
and see that $w_0 , \ldots , w_{k-1}$ is legal as well and
maximal if $v_0 , \ldots , v_{k-1}$ is.  If $i > 0$, first apply
the DSMP to $v_0 , \ldots , v_{i-1}$ and then use the same
argument.  This shows that equivalent sequences are either both
maximal popping sequences or neither.  (The case of $k=\infty$ is
trivial, since infinite popping sequences are automatically
maximal.)

To prove the lemma, we induct on $l(\omega)$, and then deal with
the case of $l(\omega)=\infty$. It is clear when $l=0$. Assuming
the lemma for $l(\omega) < L$, let $l(\omega) = L$ with maximal
popping sequence $v_0 , \ldots , v_{L-1}$.  Let $w_0 , \ldots ,
w_{k-1}$ be any other maximal popping sequence. If $w_0 = v_0$,
applying the DSMP and the induction hypothesis completes the
induction. If not, then consider the least $i$ for which $v_i =
w_0$, if any. When we pop $v_0, \ldots, v_{L-1}$, the orientation
$\{ X_{e , f(e,j)} : e \in E \}$ has a sink at $w_0$ for each $j
< i$, since the sink at $w_0$ exists until one of its edges is
popped and no other sink can contain any such edge until $w_0$ is
popped.  Thus $i$ exists and $v_j$ cannot be a neighbor of $v_i$
for $j < i$.  Hence, we can move $v_i$ to the first position by a
sequence of adjacent transpositions with non-neighbors in $G$. We
have seen that the resulting sequence $v_i , v_0 , v_1, v_2 ,
\ldots , v_{i-1} , v_{i+1} , \ldots , v_{L-1}$ is a maximal
popping sequence.  Now apply the DSMP and the induction
hypothesis to conclude that $w_0 , \ldots , w_{k-1}$ is
equivalent to $v_i , v_1 , \ldots v_{i-1} , v_{i+1} , \ldots ,
v_{L-1}$, and thus to $v_0 , \ldots , v_{L-1}$.

All that remains is the case of $l(\omega)=\infty$, i.e., the
case when all maximal popping sequences are infinite.  Given two
such sequences $v_0,v_1,\dots$ and $w_0,w_1,\dots$, the argument
{}from the previous paragraph shows that we can transform
$v_0,v_1,\dots$ so that its first element is $w_0$.  Now applying
the DSMP shows that we can bring arbitrarily long initial
segments into agreement, which is the definition of equivalence.
\end{proof}

\begin{proof}[Proof of the diamond lemma]
{}From the previous lemma, we know that $\tau = l$, so $\tau$ is
independent of the choice rule.  Furthermore, since all maximal
popping sequences are equivalent, the multisets of popped vertices
are the same.  The assertion about SFO's follows because the SFO
depends only on which vertices were popped.
\end{proof}

\begin{proof}[Proof of strong uniform time]  We prove by induction that
for any SFO $\eta$ and any finite sequence $v_0 , \ldots , v_{k-1}$,
the following event has probability $2^{-m
-\sum_{i=0}^{k-1} \deg_0(v_i)}$, where $\deg_0$ means the degree
not counting self-loops: $\tau = k$, and $v_0 , \ldots , v_{k-1}$
is a legal popping sequence for $\omega$, and $\eta (\omega) =
\eta_0$. This is vacuously true when $k = 0$. Now the probability
that the singleton $ v_0 $ is a legal pop is $2^{-\deg_0(v_0)}$,
so applying the DSMP we see that the probability of a maximal
popping sequence $v_0 , v_1 , \ldots, v_{k-1}$ with $\eta =
\eta_0$ is
$$2^{-\deg_0(v_0)} 2^{-m - \sum_{i=1}^{k-1} \deg_0(v_i)},$$
which completes the induction.

To find the probability of both $\tau = k$ and $\eta = \eta_0$
(with no restrictions on the popping sequence), we must sum this
probability over all equivalence classes of potential popping
sequences of length $k$.  We sum over equivalence classes to avoid
double counting, since for any given $\omega$,
Lemma~\ref{lem:equiv} tells us that the set of maximal popping
sequences is an equivalence class.  Since neither the summand nor
the set of sequences depends on $\eta_0$, we have proved the
lemma.
\end{proof}

\ssection{Analysis of the running time}

We still have not shown that $\tau$ is almost surely finite. While
this may appear obvious {}from some kind of Markov property, the
choice rule makes things sticky and we find it easiest to
conclude this {}from the existence of a finite upper bound on the
expected run time. To bound $\tau$ we make repeated use of the
following monotonicity principle. We let $Q(G,v)$ denote the
random number of times $v$ is popped in a maximal popping
sequence (possibly $\infty$), which, by the diamond lemma, is well
defined.

\begin{lem}[Monotonicity] \label{lem:monotonicity}
Fix $G \in \cS$ and let $H \in \cS$ be a subgraph of $G$, that is,
$V(H) \subseteq V(G)$ and $E(H) \subseteq E(G)$.  For $v \in
V(H)$,
$$
\E Q(H, v) \geq \E Q(G , v) .
$$
\end{lem}

\begin{proof}
This is proved by stochastic domination: we run sink popping
simultaneously on $H$ and $G$, using the same stacks for edges
common to both graphs. Every legal popping sequence on $G$
restricts to a legal popping sequence on $H$ as well, so under
this coupling $Q(H,v) \geq Q(G,v)$ always.
\end{proof}

\noindent{\em Remark:} Additionally, we see that equality occurs
only when no SFO on $G$ can require further popping of $v$ on $H$.

\begin{pr} \label{pr:cycles}
If $G$ is an $n$-cycle, then $\E \tau = \binom{n}{2}$.
Furthermore, conditioned on starting with $j$ edges oriented
clockwise and $n-j$ counterclockwise, the expected value of
$\tau$ is $2j(n-j)$.
\end{pr}

\begin{proof}
At any time, some of the arrows point clockwise and others
counterclockwise.  Let $Y_k$ be the number of arrows pointing
clockwise at time $k$.  Popping at any vertex causes two opposite
pointing arrows to be replaced by two random arrows.  Thus
$Y_{k+1}$ has the distribution of $Y_k + Z$ where $\Prob (Z = 1)
= \Prob (Z = -1) = 1/4$ and $\Prob (Z = 0) = 1/2$.  Therefore $\{
Y_k : k \geq 0 \}$ is a simple random walk with delay probability
of $1/2$ absorbed at $0$ and $n$. The expected absorption time
{}from $j$ is twice that for simple random walk, and thus is $2 j
(n-j)$; see equation~(3.5) on page~349 of Feller \cite{F}. Hence
$\E \tau = 2 \E Y_0 (n - Y_0)$, which is twice the expected
number of ordered pairs of edges where the first is initially
clockwise and the second initially counterclockwise.  There are
$n(n-1)$ ordered pairs of distinct edges, each having these
orientations with probability $1/4$, so $\E \tau = n(n-1)/2$.
\end{proof}

\begin{cor} \label{cor:UC}
Let $\cS_0$ denote the class of graphs in which every vertex is
in some cycle.  For $G \in \cS_0$ and $v \in V(G)$,
$$
\E Q(G,v) \leq (n-1)/2
.
$$
Equality holds for all $v$ if and only if $G$ is an $n$-cycle.
\end{cor}

\begin{proof}
Fix $G$ and $v$ and let $H$ be a cycle containing $v$.  By
monotonicity, $\E Q(G,v) \leq \E Q(H,v)$ which is at most $(n-1)
/ 2$ by Proposition~\ref{pr:cycles} and symmetry.  In fact, $\E
Q(H,v)$ is strictly less than $(n-1)/2$ unless $H$ is an
$n$-cycle.  By the remark following the proof of the monotonicity
lemma, the inequality $\E Q(G,v) \leq \E Q(H,v)$ is strict unless
no SFO on $G$ can require further popping of $v$ on $H$. In our
case, $H$ is an $n$-cycle, and $G$ is an $n$-cycle with some
chords or self-loops added.  Then (assuming $G$ is not an
$n$-cycle), there is always an SFO on $G$ that does not restrict
to an SFO on $H$: if $G$ has a self-loop, one can choose an SFO
on $G$ such that $H$ has a sink there;  if $G$ has a chord, one
can use the chord to create a short circuit across $H$ giving a
cycle of length less than $n$, orient this cycle in a loop, and
orient the other edges in $G$ towards the cycle. If $v$ is a sink
in the restriction to $H$ of such an SFO on $G$, then strict
inequality holds for $v$ (because with positive probability, sink
popping on $G$ will produce this SFO, and $v$ will still need to
be popped in $H$).
\end{proof}

\begin{lem} \label{lem:4}
For every $G \in \cS$ with $n$ vertices, and each $v \in V(G)$,
$$\E Q(G , v) \leq n-1,$$
and equality holds only when $v$ is the vertex furthest from the
self-loop in an $n$-lollipop.
\end{lem}

\begin{proof}
We induct on $G$. The base step is $G \in \cS_0$, which is
immediate {}from the previous corollary. Assume for induction
that the conclusion holds for all subgraphs of $G$. There are
three cases other than the base step.

Case 1.  $G$ is not connected.  Then the result follows {}from the
induction hypothesis and the monotonicity lemma applied to the
component $H$ of $G$ containing $v$.  Equality never occurs.

If $G$ is connected and not in $\cS_0$, then $G$ must contain an
isthmus, i.e., an edge whose removal disconnects $G$.

Case 2.  Some edge $e$ disconnects $G$ into two components both
in $\cS$. Again the result follows {}from the induction hypothesis
applied to the component $H$ of $G \setminus \{e\}$ containing
$v$, and equality never occurs.

Case 3.  $G$ has an isthmus and removal of any isthmus always
leaves a component that is a tree.  Then $G$ has a leaf $z$.  If
$v \neq z$ then the result follows immediately {}from monotonicity
with $H = G \setminus \{z\}$.  If $v$ is the only leaf, then let
$w$ be its neighbor. Choose a popping order that pops $v$ whenever
possible, and otherwise executes any choice rule for sink popping
on $H := G \setminus \{v\}$. Initially there is a $1/2$ chance
that $v$ is a sink, in which case it is popped a mean $2$
geometric number of times until the edge $vw$ points to $w$.
Then, each time $w$ is popped, the probability is $1/2$ that this
edge is reversed, in which case it takes another mean $2$
geometric number of pops to reverse it again.  Thus
$$\E Q(G,v) = 1 + \E Q(H , w) . $$
By induction, this is at most $1 + (n-2)$. Equality occurs for
$v$ in $G$ if and only if it occurs for $w$ in $H$, so we see by
induction that it holds only at the end of a lollipop.
\end{proof}

\begin{proof}[Proof of Lemma~\ref{lem:n choose 2}]
We prove the lemma by induction, following the pattern of the
last proof.  The base step is $G \in \cS_0$, in which case the
lemma follows {}from Corollary~\ref{cor:UC}. In the cases 1 and 2
of the induction, if $G$ is disconnected or the union of two
graphs in $\cS$ along an added edge, the result is again
immediate {}from the subadditivity of the function $n \mapsto
\binom{n}{2}$ and monotonicity.  Finally, if $G$ has a leaf $v$,
we set $H := G \setminus \{v\}$ and observe that the number of
pops $\tau_G$ and $\tau_H$ on $G$ and $H$ respectively are
related by $\tau_G = \tau_H + Q(G,v)$. Thus
$$\E \tau_G = \E \tau_H + \E Q(G,v) \leq \binom{n-1}{2} + (n-1)
   = \binom{n}{2} \, .$$
By the previous lemma, the last inequality is strict unless $H$
is an $n$-lollipop and its vertex of degree 1 is the neighbor of
$v$ in $G$. This completes the induction.
\end{proof}

\begin{proof}[Proof of Theorem~\ref{th:main}]
The theorem follows immediately {}from combining
Lemmas~\ref{lem:diamond}, \ref{lem:unif}, \ref{lem:n choose 2},
and~\ref{lem:4}.
\end{proof}

\ssection{Further proofs}

The $n$-cycle and $n$-lollipop have the worst mean run times.
Here we prove Proposition~\ref{pr:equal}, namely that the run
time distributions are in fact identical.

\begin{proof}[Proof of Proposition~\ref{pr:equal}]
Number the vertices of the $n$-lollipop $0 , \ldots , n-1$ with
$0$ being the leaf.  Always pop the sink with lowest number.  Let
$Y_k$ denote the sink popped at time $k$.  Clearly $Y_0$ is $-1$
plus a mean $2$ geometric random variable, with the proviso that a
value of $n-1$ or higher represents the terminal state in which
no sink needs to be popped. Let $\F_k$ be the $\sigma$-field
generated by $Y_0 , \ldots , Y_{k-1}$. We claim that $\{ Y_k : k
\geq 0 \}$ is a time-homogeneous Markov chain with respect to $\{
\F_k \}$ and that {}from any state $j > 0$ its increments are
$-2$ plus a mean $2$ geometric, jumping to the terminal state if
it reaches $n-1$ or greater, and {}from state $0$ the same thing
with $-2$ replaced by $-1$ (thus the jump {}from $0$ is resampled
if it hits $-1$). All that is needed to check this is an
inductive verification that the orientations of edges between
vertices of higher index than $Y_k$ are conditionally i.i.d.\ fair
coin flips given $\F_k$, which is straightforward.

Now we show that the running time on an $n$-cycle is also equal
to the time for a random walk to hit at least $n-1$ when its
increments are $-2$ plus a mean $2$ geometric, resampled if it
hits $-1$.  At time $k$, let $Y_k$ denote the least index of a
sink when the edge {}from $n-1$ to $0$ is oriented toward $0$ and
$n-1$ minus the greatest index of a sink when the edge is oriented
toward $n-1$. In other words, this quantity is the distance {}from
the head of the $0,n-1$ edge to the nearest sink in that
direction.  We always choose the pop that sink. The only time the
$0,n-1$ edge can change orientations is when $Y_k = 0$, in which
case $Y_{k+1}$ will be $-1$ plus a mean $2$ geometric; when $Y_k >
0$ verification of the conditional increment is trivial. The
stopping rule is, again, that one must jump to $n-1$ or greater,
and $Y_0$ has the right distribution for the same reason as
before, so the sequence has the same distribution.
\end{proof}

Our final result deals with the run time started {}from an
arbitrary state, that is, the conditional distribution of $\tau$
given $\F_0$ (as defined in Proposition~\ref{pr:extremal}). While
this quantity is a hidden variable as far as users of the
algorithm are concerned, it has relevance to the distribution of
the run time, as well as having some intrinsic interest.  We
begin again with a result on the $n$-cycle.

\begin{pr} \label{pr:alsace}
Let $G$ be an $n$-cycle.  Then for every $v \in V(G)$,
$$
\E (Q(G,v) \,|\, \F_0) \leq 3n/4,
$$
with equality if and only if $n$ is even and all edges are
oriented along the direction of shortest travel to $v$.
\end{pr}

\begin{proof}
Number the vertices $0 , \ldots , n-1$ mod $n$. We first
establish that the discrete Laplacian of $\E (Q(G,v) \,|\, \F_0)$
depends on the initial orientation of $G$ via
\begin{equation} \label{eq:2}
\E \left [ \left. Q(G,v) - \frac{Q(G,v+1) +
Q (G,v-1)}{2} \,\right|\, \F_0 \right ]
=
\begin{cases}
\phantom{-}1 & \textup{if $v$ is a sink,}\\
-1 & \textup{if $v$ is a source, and}\\
\phantom{-}0 & \textup{otherwise.}
\end{cases}
\end{equation}
To see this, choose any popping order and let $Y(v,k)$ denote the
in-degree of $v$ at time $k$, that is, the number of $e \in E(G)$
adjacent to $v$ for which $X_{e , f(e,k)}$ is oriented toward
$v$. Then, conditionally on anything up to time $k$,
$$
\E Y(v,k+1) = \E Y(v,k) - \Prob (v_k = v) +
\frac{\Prob (v_k = v+1) +
   \Prob (v_k = v-1)}{2} \, ,
$$
since any pop at $v$ reduces the expected in-degree by mean 1 and
any pop at a neighbor of $v$ increases it by mean $1/2$.  Summing
over $k$, conditioning on $\F_0$ and using $Y(v,\tau) \equiv 1$
proves~(\ref{eq:2}).

{}From Proposition~\ref{pr:cycles} we know that $\E (\tau \,|\,
\F_0) = 2Y_0 (n - Y_0)$ where $Y_0$ is the number of initial
clockwise arrows (edges oriented {}from $i+1$ to $i$ mod $n$ for
some $i$). This, along with~(\ref{eq:2}), determines $\E (Q(G,
\cdot) \,|\, \F_0)$, since the difference of any two candidates
for this function would be a harmonic function on the cycle, and
hence constant.

In general,
\begin{equation} \label{eq:formula}
\E (Q(G,v) \,|\, \F_0) =
\frac{k}{n} \left(1 + 3n - 2k - \frac{2}{k}
\sum_{j=1}^k a_j\right),
\end{equation}
if there are $k$ clockwise edges pointing {}from $v + a_j$ to $v +
a_j - 1$ for a set $\{ a_1 , \ldots , a_k \} \subseteq \{ 1 ,
\ldots n \}$ (addition taken mod $n$). To prove this formula, we
need only prove that the right hand side satisfies the two
properties that characterize the left hand side. The sum over all
$v$ (i.e., $\E(\tau \,|\, \F_0)$) is easy, since it equals
$$
k(1+3n-2k) - \frac{2}{n}\sum_{j=1}^k \sum_{i=1}^n i,
$$
which does indeed simplify to $2k(n-k)$.  To check that the right
hand side of \eqref{eq:formula} works in \eqref{eq:2}, we proceed
as follows.  Let $f(v)$ be the right hand side of
\eqref{eq:formula}. Then
$$
g(v) \deq f(v+1)-f(v) = -\frac{2}{n}
\left(-k + n \delta_{\min \{a_j\},1} \right),
$$
where $\delta$ is the Kronecker delta. Hence,
\begin{equation*}
-\frac{g(v)-g(v-1)}{2} =
\begin{cases}
\phantom{-}1 & \textup{if $v$ is a sink,}\\
-1 & \textup{if $v$ is a source, and}\\
\phantom{-}0 & \textup{otherwise,}
\end{cases}
\end{equation*}
as desired.

Equation \eqref{eq:formula} makes it easy to see when $\E (Q(G,v)
\,|\, \F_0)$ is maximized: that can occur only when
$\{a_1,\ldots,a_k\} = \{1,\ldots,k\}$, in which case $\E (Q(G,v)
\,|\, \F_0)$ equals $3n(k/n)(1-k/n)$.  This quantity is bounded
above by $3n/4$, with equality if and only if $n=2k$.
\end{proof}

\begin{proof}[Proof of Proposition~\ref{pr:extremal}]
Induct again, as in the proof of Lemma~\ref{lem:4} and the main
theorem. Simultaneously, we show by induction that $\E (Q(G,v)
\,|\, \F_0) \leq 2(n-1)$, with equality only for the leaf of an
$n$-lollipop and initial conditions $X_{e,0}$ all pointing toward
$v$.  Note that the previous proposition proves this bound for all
$G \in \cS_0$, if we use monotonicity (which also holds
conditioned on the initial orientation).

For the base case, $G \in \cS_0$ and Proposition~\ref{pr:alsace}
shows that in fact $\E (\tau \,|\, \F_0) < 3 n^2 / 4$.  There is
strict inequality because equality in Proposition~\ref{pr:alsace}
cannot hold simultaneously for all vertices in a cycle.  It
follows that $\E (\tau \,|\, \F_0) < n(n-1)$ unless $n \le 3$.
The cases with $n \le 3$ are easily dealt with: those with $n=2$
are trivial, and for $n=3$ the worst case contains a $3$-cycle,
which can be analyzed using the sharper bounds in the proof of
Proposition~\ref{pr:alsace}.

When $G$ is not connected or is the union of two graphs in $\cS$
along an added isthmus, the $n(n-1)$ bound is immediate {}from
subadditivity of $n(n-1)$ and monotonicity, and the $2(n-1)$
bound follows from monotonicity. Finally, when $G$ has a leaf
$v$, set $H := G \setminus \{v\}$ as before. This time, in the
worst case we know that $v$ is a sink initially, so $\E (Q(G , v)
\,|\, \F_0)$ is bounded by $2 + \E (Q(H , w) \,|\, \F_0)$.  This
verifies the conclusion that $\E (Q(G,v) \,|\, \F_0) \leq
2(n-1)$, and adding this to $\E \tau_H$ gives, by induction, at
most $(n-1) (n-2) + 2 (n-1) = n(n-1)$, which completes the proof
of the upper bound;  the conditions for equality are clear from
the proof.
\end{proof}

\ssection{Questions}

It is tempting to view both cycle popping and sink popping as
special cases of what might be called ``partial rejection
sampling:'' to generate a random structure, choose a random
candidate, and if it has any flaws, locally rerandomize until it
is flawless. Does partial rejection sampling apply to other
natural combinatorial problems? Can one develop a general theory?
Note that Fill and Huber's randomness recycler \cite{FH} also uses
the idea of rejecting only part of a structure, although in a
different way.

Cycle popping was applied to the study of random spanning trees
on $\Z^d$, as well as some more general graphs, in \cite{BLPS}.
It would be interesting if sink popping could be used similarly.
Do random sink-free orientations on $\Z^d$ exhibit any interesting
or surprising structure?

\ssections{Acknowledgements}

We thank David Wilson and Peter Winkler for helpful discussions,
and the anonymous referee for useful comments on the manuscript.


\begin{thebibliography}{BLPS}

\vskip 2.8ex

\bibitem[A]{A} D.~Aldous, {\sl A random walk construction of
uniform spanning trees and uniform labelled trees\/}, SIAM
Journal of Discrete Mathematics {\bf 3 (4)} (1990), 450--465.

\bibitem[BLPS]{BLPS} I.~Benjamini, R.~Lyons, Y.~Peres, and
O.~Schramm, {\sl Uniform spanning forests\/}, Ann.\ Probab.\ {\bf
29} (2001), no.\ 1, 1--65.

\bibitem[B]{B} A.~Broder, {\sl Generating random spanning trees\/},
Proceedings of 30th Annual Symposium on Foundations of Computer
Science (Research Triangle Park, NC, 1989), 442--447, IEEE, 1989.

\bibitem[BD]{BD}
R.~Bubley and M.~Dyer, {\sl Graph orientations with no sink and an
approximation for a hard case of \#SAT\/}, Proceedings of the
Eighth Annual ACM-SIAM Symposium on Discrete Algorithms (New
Orleans, LA, 1997), 248--257, ACM, New York, 1997.

\bibitem[CEP]{CEP}
H.~Cohn, N.~Elkies, and J.~Propp, {\sl Local statistics for
random domino tilings of the Aztec diamond\/}, Duke Math.\ J.\
{\bf 85} (1996), no.\ 1, 117--166, arXiv:math.CO/0008243.

\bibitem[E]{E}
K.~Eriksson, {\sl Strong convergence and the polygon property of
$1$-player games\/}, Discrete Math.\ {\bf 153} (1996), no. 1--3,
105--122.

\bibitem[F]{F}
W.~Feller, {\sl An introduction to probability theory and its
applications\/}, Volume I, third edition, John Wiley \& Sons,
Inc., New York, 1968.

\bibitem[FH]{FH}
J.~Fill and M.~Huber, {\sl The randomness recycler: a new
technique for perfect sampling\/}, Proceedings of the Forty-first
Annual Symposium on Foundations of Computer Science (Redondo
Beach, CA, 2000), 503--511, IEEE, 2000.

\bibitem[H]{H}
M.~Huber, {\sl Exact sampling and approximate counting
techniques\/}, Proceedings of the Thirtieth Annual ACM Symposium
on the Theory of Computing (Dallas, TX, 1998), 31--40, ACM, New
York, 1998.

\bibitem[JS]{JS}
M.~Jerrum and A.~Sinclair, {\sl Approximate counting, uniform
generation and rapidly mixing Markov chains\/}, Inform.\ and
Comput.\ {\bf 82} (1989), no. 1, 93--133.

\bibitem[P]{P} R.~Pemantle,
{\sl Choosing a spanning tree for the integer lattice
uniformly\/}, Ann.\ Probab.\ {\bf 19} (1991), no. 4, 1559--1574.

\bibitem[PW1]{PW} J.~Propp and D.~Wilson,
{\sl Exact sampling with coupled Markov chains and applications
to statistical mechanics\/}, Random Structures and Algorithms
{\bf 9} (1996), no. 1--2, 223--252.

\bibitem[PW2]{PW2} J.~Propp and D.~Wilson,
{\sl How to get a perfectly random sample {}from a generic Markov
chain and generate a random spanning tree of a directed graph\/},
Journal of Algorithms {\bf 27} (1998), 170--217.

\bibitem[W]{W}
D.~Wilson, {\sl Generating random spanning trees more quickly
than the cover time\/}, Proceedings of the Twenty-eighth Annual
ACM Symposium on the Theory of Computing (Philadelphia, PA,
1996), 296--303, ACM, New York, 1996.

\end{thebibliography}
\end{document}